\documentclass[a4paper,12pt]{amsart}
\usepackage{amsfonts}
\usepackage{amssymb}
\usepackage{ifthen}
\usepackage{amscd}
\usepackage{amsxtra}
\usepackage{graphicx}
\usepackage{color}
\nonstopmode \numberwithin{equation}{section}
\newtheorem{thm}{Theorem}
\newtheorem{lem}{Lemma}
\newtheorem{cor}{Corollary}[section]

\newtheorem{cl}{Claim}
\newtheorem{ca}{Case}
\newtheorem{sca}{Subcase}
\newtheorem{scl}{Subclaim}
\newtheorem{conj}{Conjecture}

\theoremstyle{definition}
\newtheorem{defn}{Definition}

\newtheorem{op}[equation]{Open Problem}
\newtheorem{ques}[equation]{Question}
\newtheorem{rem}{Remark}[section]
\newtheorem{exam}[equation]{Example}

\newcounter {own}
\def\theown {\thesection       .\arabic{own}}

\newenvironment{pf}[1][]{%
 \vskip 3mm
 \noindent
 \ifthenelse{\equal{#1}{}}%
  {{\slshape Proof. }}%
  {{\slshape #1.} }%
 }%
{\qed\bigskip}

\newcounter{alphabet}
\newcounter{tmp}
\newenvironment{Thm}[1][]{\refstepcounter{alphabet}%
\bigskip%
\noindent%
{\bf Theorem \Alph{alphabet}}%
\ifthenelse{\equal{#1}{}}{}{ (#1)}%
{\bf .} \itshape}{\vskip 8pt}

\makeatletter
\newcommand{\Ref}[1]{\@ifundefined{r@#1}{}{\setcounter{tmp}{\ref{#1}}\Alph{tmp}}}
\makeatother

\newenvironment{Lem}[1][]{\refstepcounter{alphabet}%
\bigskip%
\noindent%
{\bf Lemma \Alph{alphabet}}%
{\bf .} \itshape}{\vskip 8pt}

\newcommand{\ID}{{\mathbb D}}

\def\be{\begin{equation}}
\def\ee{\end{equation}}

\newcommand{\bee}{\begin{enumerate}}
\newcommand{\eee}{\end{enumerate}}

\newcommand{\blem}{\begin{lem}}
\newcommand{\elem}{\end{lem}}
\newcommand{\bthm}{\begin{thm}}
\newcommand{\ethm}{\end{thm}}
\newcommand{\bcor}{\begin{cor}}
\newcommand{\ecor}{\end{cor}}
\newcommand{\beg}{\begin{exam}}
\newcommand{\eeg}{\end{exam}}
\newcommand{\begs}{\begin{examples}}
\newcommand{\eegs}{\end{examples}}
\newcommand{\bdefe}{\begin{defn}}
\newcommand{\edefe}{\end{defn}}
\newcommand{\bprob}{\begin{prob}}
\newcommand{\eprob}{\end{prob}}
\newcommand{\bques}{\begin{ques}}
\newcommand{\eques}{\end{ques}}
\newcommand{\bei}{\begin{itemize}}
\newcommand{\eei}{\end{itemize}}
\newcommand{\bcon}{\begin{conj}}
\newcommand{\econ}{\end{conj}}
\newcommand{\bop}{\begin{op}}
\newcommand{\eop}{\end{op}}

\newcommand{\bca}{\begin{ca}}
\newcommand{\eca}{\end{ca}}
\newcommand{\bsca}{\begin{sca}}
\newcommand{\esca}{\end{sca}}

\newcommand{\bcl}{\begin{cl}}
\newcommand{\ecl}{\end{cl}}

\newcommand{\bscl}{\begin{scl}}
\newcommand{\escl}{\end{scl}}

\newcommand{\bcons}{\begin{conjs}}
\newcommand{\econs}{\end{conjs}}
\newcommand{\bprop}{\begin{propo}}
\newcommand{\eprop}{\end{propo}}
\newcommand{\br}{\begin{rem}}
\newcommand{\er}{\end{rem}}
\newcommand{\brs}{\begin{rems}}
\newcommand{\ers}{\end{rems}}
\newcommand{\bo}{\begin{obser}}
\newcommand{\eo}{\end{obser}}
\newcommand{\bos}{\begin{obsers}}
\newcommand{\eos}{\end{obsers}}
\newcommand{\bpf}{\begin{pf}}
\newcommand{\epf}{\end{pf}}
\newcommand{\ba}{\begin{array}}
\newcommand{\ea}{\end{array}}
\newcommand{\beq}{\begin{eqnarray}}
\newcommand{\beqq}{\begin{eqnarray*}}
\newcommand{\eeq}{\end{eqnarray}}
\newcommand{\eeqq}{\end{eqnarray*}}

\newcounter{minutes}
\divide\time by 60
\newcounter{hours}
\multiply\time by 60 \addtocounter{minutes}{-\time}

\begin{document}

\bibliographystyle{amsplain}
\title []
{Lengths, areas  and Lipschitz-type spaces of planar   harmonic
mappings}

\def\thefootnote{}
\footnotetext{ \texttt{\tiny File:~\jobname .tex,
          printed: \number\day-\number\month-\number\year,
          \thehours.\ifnum\theminutes<10{0}\fi\theminutes}
} \makeatletter\def\thefootnote{\@arabic\c@footnote}\makeatother

\author{Sh. Chen}
\address{Sh. Chen, Department of Mathematics and Computational
Science, Hengyang Normal University, Hengyang, Hunan 421008,
People's Republic of China.} \email{mathechen@126.com}

\author{S. Ponnusamy $^\dagger $ 
}
\address{S. Ponnusamy,
Indian Statistical Institute (ISI), Chennai Centre, SETS (Society
for Electronic Transactions and security), MGR Knowledge City, CIT
Campus, Taramani, Chennai 600 113, India. }
\email{samy@isichennai.res.in, samy@iitm.ac.in}

\author{ A. Rasila }
\address{A. Rasila, Department of Mathematics and Systems Analysis, Aalto University, P. O. Box 11100, FI-00076 Aalto,
 Finland.} \email{antti.rasila@iki.fi}


\keywords{Harmonic mapping, three circles theorem,  area function.\\
$
^\dagger${\tt ~~Corresponding author. This author is on leave from the Department of Mathematics,
Indian Institute of Technology Madras, Chennai-600 036, India}
}

\begin{abstract}
In this paper, we establish a three circles type
theorem, involving the harmonic area function, for harmonic mappings. Also, we give bounds for
length and area distortion for harmonic quasiconformal mappings. Finally, we will study certain
Lipschitz-type spaces on harmonic mappings.
\end{abstract}


\maketitle \pagestyle{myheadings} \markboth{ Sh. Chen,  S. Ponnusamy
and A. Rasila }{Lengths, areas  and Lipschitz-type spaces of planar
harmonic mappings}

\section{Introduction and main results }\label{csw-sec1}
Let $D$ be a simply connected subdomain of the complex plane
$\mathbb{C}$. A complex-valued function $f$ defined in $D$ is called
a {\it harmonic mapping} in $D$ if and only if both the real and the
imaginary parts of $f$ are real harmonic in $D$. It is known that
every harmonic mapping $f$ defined in $D$ admits a decomposition
$f=h+\overline{g}$, where $h$ and $g$ are analytic in $D$. Since the
Jacobian $J_f$ of $f$ is given by
$$J_f=|f_{z}|^2-|f_{\overline{z}}|^2:=|h'|^2-|g'|^2,
$$
$f$ is locally univalent and sense-preserving in  $D$ if and only if
$|g'(z)|<|h'(z)|$ in $D$; or equivalently if $h'(z)\neq0$ and the
dilatation $\omega =g'/h'$ has the property that $|\omega(z)|<1$ in $D$ (see
\cite{Lewy}). Let ${\mathcal H}(D)$ denote the class of all sense-preserving
harmonic mappings in $D$.
We refer to \cite{Clunie-Small-84,Du} for basic results in the theory
of planar harmonic mappings.

For $a\in\mathbb{C}$, let $\ID(a,r)=\{z:\, |z-a|<r\}$. In
particular, we use $\mathbb{D}_r$ to denote the disk
$\mathbb{D}(0,r)$ and  $\mathbb{D}$, the open unit disk $\ID_1$.


The classical theorem of three circles \cite{A,R}, also called {\it
Hadamard's three circles theorem}, states that if $f$ is an analytic
function in the annulus
$B(r_{1},r_{2})=\{z:\, 0<r_{1}<|z|=r<r_{2}<\infty\}$,
continuous on $\overline{B(r_{1},r_{2})}$, and $M_{1}$, $M_{2}$
and $M$ are the maxima of $f$ on the three circles corresponding to
$r_{1}$, $r_{2}$ and $r$, respectively, then
$$M^{\log \frac{r_{2}}{r_{1}}}\leq M_{1}^{\log \frac{r_{2}}{r}}M_{2}^{\log \frac{r}{r_{1}}}.
$$
Equivalently, we can reformulate this result into a simpler form.
That is if $f$ is analytic on the annulus
$B(r_{1},1)=\{z:\, 0<r_{1}<|z|<1\}$, continuous on the closure, and
$$|f(z)|\leq m=r_{1}^{\alpha},~|z|=r_{1}~\mbox{ and }~|f(z)|\leq1,~|z|=1,
$$
then Hadamard's result states that, for $r_{1}\leq r\leq1,$
$$|f(z)|\leq m^{\frac{\log r}{\log r_{1}}}=r^{\alpha},~|z|=r,
$$
where $\alpha$ is an integer.

The original three circles theorem was given by Hadamard without proof in 1896 \cite{Hada},
 and comprehensive discussion about the history of this result can be found  in
 \cite[pp. 323--325]{Maz} and \cite{R}. It is a natural question, what results of this type can be proved for other classes of functions
 and, indeed, there are numerous generalizations of the thee circles theorem in the literature, see e.g.
\cite{Brummel,MRV,Protter,Vyborny}. In this paper, our first aim is to establish an area version of the
three circles theorem (cf. area version of Schwarz' lemma \cite{Burckel}).

For a harmonic mapping $f$ in $\mathbb{D}$ and $r\in[0,1)$, the {\it
harmonic area function} $S_{f}(r)$ of $f$, counting multiplicity, is
defined by
$$S_{f}(r)=\int_{\mathbb{D}_{r}}J_f(z)\,d\sigma(z),
$$
where  $d\sigma$ denotes the normalized Lebesgue area measure on $\mathbb{D}$ (cf.
\cite{CPR}). In particular, let
$$S_{f}(1)=\sup_{0<r<1}S_{f}(r).
$$

\begin{thm}\label{thm-1r}
Let $f=h+\overline{g}$ be harmonic in $\mathbb{D}$, where $h$ and
$g$ are analytic. If $S_{f}(r_{1})\leq m<1$, $S_{f}(1)\leq 1$ and
for all $n\in\{1,2,\ldots\}$, $|g^{(n)}(0)|\leq|h^{(n)}(0)|$,  then
for $r_{1}\leq r<1,$
\be\label{eq-1} S_{f}(r)\leq m^{\frac{\log r}{\log r_{1}}}.
\ee
The estimate of \eqref{eq-1} is sharp and  the extremal function is
$f(z)= \alpha z+\beta\overline{z}$, where $\alpha$ and $\beta$ are constant with
$|\alpha|^{2}-|\beta|^{2}=1$.
\end{thm}

\begin{cor} Let $f$ be analytic in $\mathbb{D}$
satisfying $S_{f}(r_{1})\leq m$ and $S_{f}(1)\leq 1,$ where $0<
r_{1}<1$. Then for $r_{1}\leq r<1,$
\be\label{eq-2}
S_{f}(r)\leq m^{\frac{\log r}{\log r_{1}}}.
\ee
The estimate of \eqref{eq-2} is sharp and  the extremal function is $f(z)=\lambda z$, where
$|\lambda|=1$ are constant.
\end{cor}


For $p\in(0,\infty]$, the {\it harmonic Hardy space $h^{p}$}
consists of all harmonic functions $f$ such that  $
\|f\|_{p}<\infty$, where
$$\|f\|_{p}= \begin{cases}
\displaystyle\sup_{0<r<1}M_{p}(r,f)
& \mbox{ if } p\in(0,\infty),\\
\displaystyle\sup_{z\in\mathbb{D}}|f(z)| &\mbox{ if } p=\infty,
\end{cases}
~\mbox{ and }~
M_{p}^{p}(r,f)=\frac{1}{2\pi}\int_{0}^{2\pi}|f(re^{i\theta})|^{p}\,d\theta.
$$
If $f\in h^{p}$ for some $p>0$, then the radial limits
$$f(e^{i\theta})=\lim_{r\rightarrow1-}f(re^{i\theta})
$$
exist for almost every $\theta\in[0,2\pi)$ (cf. \cite{Du}).

We recall that a function $f\in {\mathcal H}(\ID)$ is said to be {\it
$K$-quasiregular}, $K\in[1,\infty)$, if for $z\in\ID$,
$\Lambda_{f}(z)\leq K\lambda_{f}(z)$.
In addition, if $f$ is univalent in $\ID$, then $f$ is called a
{\it $K$-quasiconformal} harmonic mapping $\ID$.

Let $\Omega$ be a domain of $\mathbb{C}$, with non-empty boundary.
Let $d_{\Omega}(z)$ be the Euclidean distance from $z$ to the
boundary $\partial \Omega$ of $\Omega$. In particular, we always use
$d(z)$ to denote the Euclidean distance from $z$ to the boundary of
$\mathbb{D}.$ The area of a set $G\subset\mathbb{C}$ is denoted by
$A(G)$. The area problem of analytic functions has attracted
much attention (see \cite{ATU, U,Y1,Y}). We investigate the area
problem of harmonic mappings and obtain the following result.

\begin{thm}\label{thm-2}
Let $\Omega_1$ and $\Omega _2$ be two proper and simply connected
subdomains of $\mathbb{C}$ containing  the point of origin. Then for
a sense-preserving and $K$-quasiconformal harmonic mapping $f$
defined in $\Omega_{1}$ with $f(0)=0$, \be\label{eq-t1} K
A\big(f(\Omega_{1})\cap\Omega_{2}\big)+A(f^{-1}(\Omega_{2}))\geq
\min\{d_{\Omega_{1}}^{2}(0),d_{\Omega_{2}}^{2}(0)\}. \ee Moreover,
if $K=1$, then the estimate of \eqref{eq-t1} is sharp.
\end{thm}

We remark that Theorem \ref{thm-2} is a generalization of \cite[Theorem]{U}.

 For a harmonic mapping $f$ defined on $\mathbb{D}$, we use the following
standard notations:
$$\Lambda_{f}(z)=\max_{0\leq \theta\leq 2\pi}|f_{z}(z)+e^{-2i\theta}f_{\overline{z}}(z)|
=|f_{z}(z)|+|f_{\overline{z}}(z)|
$$
and
$$\lambda_{f}(z)=\min_{0\leq \theta\leq 2\pi}|f_{z}(z)+e^{-2i\theta}f_{\overline{z}}(z)|
=\big | \, |f_{z}(z)|-|f_{\overline{z}}(z)|\, \big |.
$$
Further, a planar harmonic mapping $f$ defined on $\mathbb{D}$ is called a {\it
harmonic Bloch mapping} if
$$\beta_{f}=\sup_{ z,w\in\mathbb{D},\ z\neq w}\frac{|f(z)-f(w)|}{\rho(z,w)}<\infty.
$$
Here $\beta_{f}$ is called the {\it Lipschitz number} of $f$, and
$$\rho(z,w)=\frac{1}{2}\log\left(\frac{1+|(z-w)/(1-\overline{z}w)|}
{1-|(z-w)/(1-\overline{z}w)|}\right)=\mbox{arctanh}\left |\frac{z-w}{1-\overline{z}w}\right |
$$
denotes the hyperbolic distance between $z$ and $w$ in $\mathbb{D}$.
It is known that
$$\beta_{f}=\sup_{z\in\mathbb{D}}\big\{(1-|z|^{2})\Lambda_{f}(z)\big\}.
$$
Clearly, a harmonic Bloch mapping $f$ is uniformly continuous as a
map between metric spaces,
$$f:\,(\mathbb{D}, \rho) \to ( \mathbb{C}, | \cdot|),
$$
and for all $z,w \in \mathbb{D}$ we have the Lipschitz inequality
$$|f(z)-f(w)| \le \beta_{f} \, \rho(z,w) \,.
$$
A well-known fact is that the set of all harmonic Bloch mappings, denoted by the
symbol $\mathcal{HB}$, forms a complex Banach space with the norm
$\|\cdot\|$ given by
$$\|f\|_{\mathcal{HB}}=|f(0)|+\sup_{z\in\mathbb{D}}\{(1-|z|^{2})\Lambda_{f}(z)\}.
$$
Specially, we use  $\mathcal{B}$  to denote the set of all analytic
functions defined in $\mathbb{D}$ which forms a complex Banach space
with the norm
$$\|f\|_{\mathcal{B}}=|f(0)|+\sup_{z\in\mathbb{D}}\{(1-|z|^{2})|f'(z)|\}.
$$
The reader is referred to \cite[Theorem 2]{Co} (or \cite{CPW0,CPW1})
for a detailed discussion.

For $r\in[0,1)$,   the length of the curve
$C(r)=\big\{w=f(re^{i\theta}):\, \theta\in[0,2\pi]\big\}$, counting
multiplicity, is defined by
$$l_{f}(r)=\int_{0}^{2\pi}|df(re^{i\theta})|= r\int_{0}^{2\pi}\left|f_{z}(re^{i\theta})-e^{-2i\theta}
f_{\overline{z}}(re^{i\theta})\right|d\theta,
$$
where  $f$ is a harmonic mapping defined in $\mathbb{D}$. In particular, let
$l_{f}(1)=\sup_{0<r<1}l_{f}(r)$.

\begin{thm}\label{thm-4}
Let $f(z)=\sum_{n=0}^{\infty}a_{n}z^{n}+\sum_{n=1}^{\infty}\overline{b}_{n}\overline{z}^{n}$
be a sense-preserving $K$-quasiconformal harmonic mapping. If
$l_{f}(1)<\infty$, then for $n\geq1$,
\be\label{thm-eq1}
|a_{n}|+|b_{n}|\leq\frac{Kl_{f}(1)}{2n\pi}
\ee
and
\be\label{thm-eq2}
\Lambda_{f}(z)\leq
\frac{l_{f}(1)\sqrt{K}}{2\pi(1-|z|)}.
\ee
Moreover, $f\in\mathcal{HB}$ and $\beta_{f}\leq\frac{l_{f}(1)\sqrt{K}}{\pi}.$ In
particular, if $K=1$, the estimates of \eqref{thm-eq1} and
\eqref{thm-eq2} are sharp, and the extremal function is
$f(z)=z.$
\end{thm}


A continuous increasing function $\omega:\, [0,\infty)\rightarrow
[0,\infty)$ with $\omega(0)=0$ is called a {\it majorant} if
$\omega(t)/t$ is non-increasing for $t>0$. Given a subset $\Omega$
of $\mathbb{C}$, a function $f:\, \Omega\rightarrow \mathbb{C}$ is
said to belong to the {\it Lipschitz space $L_{\omega}(\Omega)$} if
there is a positive constant $C$ such that
\be\label{eq1}
|f(z)-f(w)|\leq C\omega(|z-w|) ~\mbox{ for all $z,\ w\in\Omega.$}
\ee
For $\delta_{0}>0$, let
\be\label{eq2}
\int_{0}^{\delta}\frac{\omega(t)}{t}\,dt\leq C\cdot\omega(\delta),\
0<\delta<\delta_{0},
\ee
and
\be\label{eq3}
\delta\int_{\delta}^{+\infty}\frac{\omega(t)}{t^{2}}\,dt\leq
C\cdot\omega(\delta),\ 0<\delta<\delta_{0},
\ee
where $\omega$ is a majorant and $C$ is a positive constant.

A majorant $\omega$ is said to be {\it regular} if it satisfies the
conditions (\ref{eq2}) and (\ref{eq3}) (see
\cite{D,D1,P,Pav1,Pav2}).

Let $G$ be a proper subdomain of $\mathbb{C}$. We say that a
function $f$ belongs to the {\it local Lipschitz space }
$\mbox{loc}L_{\omega}(G)$ if (\ref{eq1}) holds, with a fixed
positive constant $C$, whenever $z\in G$ and
$|z-w|<\frac{1}{2}d_{G}(z)$ (cf. \cite{GM,L}). Moreover, $G$ is said
to be a {\it $L_{\omega}$-extension domain} if
$L_{\omega}(G)=\mbox{loc}L_{\omega}(G).$ The geometric
characterization of $L_{\omega}$-extension domains was first given
by Gehring and Martio \cite{GM}. Then Lappalainen \cite{L} generalized their characterization, and proved that $G$ is a
$L_{\omega}$-extension domain if and only if each pair of points
$z,w\in G$ can be joined by a rectifiable curve $\gamma\subset G$
satisfying
\be\label{eq1.0}
\int_{\gamma}\frac{\omega(d_{G}(z))}{d_{G}(z)}\,ds(z) \leq C\omega(|z-w|)
\ee
with some fixed positive constant $C=C(G,\omega)$, where $ds$ stands for the
arc length measure on $\gamma$.  Furthermore, Lappalainen \cite[Theorem 4.12]{L}
proved that $L_{\omega}$-extension domains  exist only for majorants
$\omega$ satisfying  (\ref{eq2}).

\begin{Thm}{\rm (\cite[Theorem 3]{Ho})}\label{ThmA2}
$f\in\mathcal{B}$ if and only if
$$\sup_{z,w\in\mathbb{D},z\neq w}\left\{
\frac{\sqrt{(1-|z|^{2})(1-|w|^{2})}|f(z)-f(w)|}{|z-w|}\right\}<\infty.
$$
\end{Thm}

 The following result is a generalization of
Theorem \Ref{ThmA2}. For the related studies of this topic for real
functions, we refer to \cite{Pav,Re}.

\begin{thm}\label{thm-CPW}
Let $f$ be a harmonic mapping in $\mathbb{D}$ and  $\omega$ be a
majorant. Then the following are equivalent:

\item{{{\rm(a)}}}~There exists a constant $C_1>0$ such that for all $z\in\mathbb{D}$,
$$\Lambda_{f}(z)\leq C_1\omega\left(\frac{1}{d(z)}\right);
$$
\item{{{\rm(b)}}}~There exists a constant $C_2>0$ such that for all $z, w\in\mathbb{D}$ with $z\neq w$,
$$\frac{|f(z)-f(w)|}{|z-w|}\leq
C_2\omega\left(\frac{1}{\sqrt{d(z)d(w)}}\right);
$$

\item{{{\rm(c)}}}~There exists a constant $C_3>0$ such that for all  $r\in(0,d(z)]$,
$$\frac{1}{|\mathbb{D}(z,r)|}\int_{\mathbb{B}^{n}(z,r)}|f(\zeta)-f(z)|\,dA(\zeta)\leq C_3r\omega\Big(\frac{1}{r}\Big),
$$
where $dA$ denotes the Lebesgue area measure in $\mathbb{D}$.

\end{thm}


Note that if  $\omega(t)=t$ and $f$ is analytic, then
(a)$\Longleftrightarrow$(b) in  Theorem \ref{thm-CPW}  implies that
Theorem \Ref{ThmA2}.

Krantz \cite{Kr} proved the following Hardy-Littlewood-type theorem for harmonic
functions  with respect to the majorant
$\omega(t)=\omega_{\alpha}(t)=t^{\alpha}~(0<\alpha\leq1)$.

\begin{Thm}{\rm (\cite[Theorem 15.8]{Kr})}\label{ThmA3}
Let $u$ be a real harmonic function in $\mathbb{D}$ and
$0<\alpha\leq1$. Then $u$ satisfies
$$|\nabla u(z)|\leq C\frac{\omega_{\alpha}\big(d(z)\big)}{d(z)}
~\mbox{ for all }~z\in\mathbb{D}
$$
if and only if
$$|u(z)-u(w)|\leq C\omega_{\alpha}(|z-w|)~\mbox{ for all }~z,w\in\mathbb{D}.
$$
\end{Thm}

We generalize Theorem \Ref{ThmA3} to the following form.

\begin{thm}\label{thm-1}
Let $\omega$ be a majorant satisfying  \eqref{eq2}, $\Omega$ be a
$L_{\omega}$-extension domain and $f$ be a harmonic mapping in
$\Omega$. Then there exists a constant $C_4>0$ such that
$$\Lambda_{f}(z)\leq C_4\frac{\omega\big(d_{\Omega}(z)\big)}{d_{\Omega}(z)}
~\mbox{ for all }~z\in\Omega
$$
if and only if, for some $C_5>0$,
$$|f(z)-f(w)|\leq C_5\omega(|z-w|)~\mbox{ for all $z,w\in\Omega$}.
$$
\end{thm}

The proofs of Theorems \ref{thm-1r}, \ref{thm-2} and \ref{thm-4}
will be given in Section \ref{csw-sec2}. We will show Theorems
\ref{thm-CPW} and \ref{thm-1}  in the last part of this paper.

\section{Length and area of harmonic mappings}\label{csw-sec2}

\subsection*{Proof of Theorem \ref{thm-1r}}
 Let $f=h+\overline{g}$ be harmonic in $\mathbb{D}$ with the following
expansion
$$f(z)=\sum_{n=0}^{\infty}a_{n}z^{n}+\sum_{n=0}^{\infty}\overline{b}_{n}\overline{z}^{n},
$$
where $a_{n}=\frac{h^{(n)}(0)}{n!}$ and $b_{n}=\frac{g^{(n)}(0)}{n!}$. By hypothesis, $|a_n|\geq |b_n|$.
Then
$$S_{f}(r)=\int_{\mathbb{D}_{r}}J_{f}(\zeta)\,d\sigma(\zeta)=\int_{\mathbb{D}_{r}}
\left(|h'(\zeta)|^{2}-|g'(\zeta)|^{2}\right)d\sigma(\zeta)=
\sum_{n=1}^{\infty}n\left(|a_{n}|^{2}-|b_{n}|^{2}\right)r^{2n}.
$$
For $z\in\mathbb{D}$, let $F(z)=\sum_{n=1}^{\infty}A_{n}z^{2n},$
where $A_{n}=n\left(|a_{n}|^{2}-|b_{n}|^{2}\right)$. Since
$A_{n}\geq0$, we see that the maximum of $F$ on
$\partial\mathbb{D}_{r}$ is obtained on the real axis, that is
$$S_{f}(r)=F(r)=\max_{|z|=r}|F(z)|,
$$
where $r_{1}\leq r<1$. Hence the result follows from Hadamard's theorem.

Now we are ready to prove the sharpness part. It is not difficult to see
that $S_{f}(r_{1})=r_{1}^{2}$ and $S_{f}(1)=1$, where $f(z)= \alpha
z+\beta\overline{z}$ with $|\alpha|^{2}-|\beta|^{2}=1$. Then for
$0<r_{1}\leq r<1$, $S_{f}(r)=r^{2}=m^{\frac{\log r}{\log r_{1}}}$,
where $m=r_{1}^{2}$. The proof of the theorem is complete. \qed

\begin{lem}\label{lem-1}
Let $f$ be a sense-preserving and $K$-quasiconformal harmonic
mapping in $\mathbb{D}$ with $f(0)=0.$ If $A(f(\mathbb{D}))<\infty$,
then $f\in h^{2}$ and
\be\label{eq-t3}
\|f\|_{2}^{2}\leq KA(f(\mathbb{D})).
\ee
Moreover, if $K=1$, then the estimate of \eqref{eq-t3} is sharp and the extremal function is $f(z)=z$.
\end{lem}

\bpf Let $f$ be a sense-preserving and $K$-quasiconformal harmonic
mapping in $\mathbb{D}$ with $f(0)=0.$ Then, by the definition of $S_{f}(r)$, we see that
\begin{eqnarray*}
S_{f}(1)=A(f(\mathbb{D}))&=&\int_{\mathbb{D}}J_{f}(z)\,d\sigma(z)
=\int_{\mathbb{D}}\Lambda_{f}(z)\lambda_{f}(z)\,d\sigma(z)\\
&\geq&\frac{1}{K}\int_{\mathbb{D}}\Lambda_{f}^{2}(z)\,d\sigma(z)\\
&\geq&\frac{1}{K}\int_{\mathbb{D}}\left(|f_{z}(z)|^{2}+|f_{\overline{z}}(z)|^{2}\right)d\sigma(z)\\
&=&\frac{1}{K}\sum_{n=1}^{\infty}n(|a_{n}|^{2}+|b_{n}|^{2})\\
&\geq&\frac{1}{K}\sum_{n=1}^{\infty}(|a_{n}|^{2}+|b_{n}|^{2})\\
&=&\frac{1}{K}\|f\|_{2}^{2}.
\end{eqnarray*}
Then $\|f\|_{2}^{2}\leq KA(f(\mathbb{D})).$ Furthermore, if $K=1$,
then the function $f(z)=z$ shows the estimate of \eqref{eq-t3}
is sharp. The proof of this lemma is complete. \epf

\subsection*{Proof of Theorem \ref{thm-2}}
Let $D_{1}=f(\Omega_{1})\cap\Omega_{2}$ and
$D_{2}=f^{-1}(\Omega_{2}).$ It is not difficult to see that
$D_{2}=f^{-1}(D_{1})$. Without loss of generality, we assume that
for $k=1,2$, $A(D_{k})<\infty.$ Let $E$ be the component of the open
set $D_{2}$  containing the original point. Then there is a
universal covering mapping $\varphi$ such that
 $\varphi:\,\mathbb{D}\rightarrow E$
with $\varphi(0)=0.$ Let $F=f\circ\varphi$. It is easy to see that
$F$ is also a $K$-quasiconformal harmonic mapping.  By using Lemma
\ref{lem-1}, we have
\be\label{eqt4}
\|\varphi\|_{2}^{2}\leq A(E)\leq A(D_{2})
\ee
and
\be\label{eqt5}
\|F\|_{2}^{2}\leq A(f(E))\leq K A(f(D_{2}))=K A(D_{1}),
\ee
which imply that
\be\label{eqt6}
\|\varphi\|_{2}^{2}+\|F\|_{2}^{2}\leq A(D_{2})+K A(D_{1}).
\ee
Since  $\varphi$ and $F$ belong to $h^{2}$, we conclude that the linear measure $m(\gamma)$ of
$$\gamma=\big\{\xi\in\partial\mathbb{D}:~\mbox{ both $|\varphi(\xi)|$ and $|F(\xi)|$ are finite}\big\}
$$
is $2\pi$. Let
$$\gamma_{\varphi}=\{\xi\in\gamma:\, |\varphi(\xi)|\geq d_{\Omega_{1}}(0)\}
$$
and, similarly,
$$\gamma_{F}=\{\xi\in\gamma:\, |F(\xi)|\geq d_{\Omega_{2}}(0)\}.
$$
Then
\be\label{eqt7}
\|\varphi\|_{2}^{2}\geq\frac{1}{2\pi}\int_{\gamma_{\varphi}}d_{\Omega_{1}}^{2}(0)|\,d\xi|\geq
\frac{d_{\Omega_{1}}^{2}(0)m(\gamma_{\varphi})}{2\pi}
\ee
and
\be\label{eqt8}
\|F\|_{2}^{2}\geq\frac{1}{2\pi}\int_{\gamma_{F}}d_{\Omega_{2}}^{2}(0)|\,d\xi|\geq
\frac{d_{\Omega_{2}}^{2}(0)m(\gamma_{F})}{2\pi}.
\ee

{\bf Claim.} $\gamma=\gamma_{\varphi}\cup \gamma_{F}.$

Suppose $\gamma\neq\gamma_{\varphi}\cup \gamma_{F}.$ Then there is a
$\xi_{0}\in\gamma\backslash(\gamma_{\varphi}\cup \gamma_{F})$ such
that $\varphi(\xi)\in\Omega_{1}$ with
$|\varphi(\xi)|<d_{\Omega_{1}}(0)$, and $F(\xi)\in\Omega_{2}$ with
$|F(\xi)|<d_{\Omega_{2}}(0).$ Since $f$ is continuous at
$\varphi(\xi_{0}),$  we know that
$$F(\xi)=\lim_{r\rightarrow1-}F(r\xi)=\lim_{r\rightarrow1-}f(\varphi(r\xi))=f(\varphi(\xi)).
$$
On the other hand, $l=\{\varphi(r\xi):\, r\in[0,1]\}$ is a curve
joining $0$ and $\varphi(\xi)$ in $\Omega_{1}$. Hence
$\varphi(\xi)\in E$ which is a contradiction with the covering
property of $\mathbb{D}$ induced by $\varphi$ over $E$.

Hence by (\ref{eqt7}), (\ref{eqt8})  and the Claim, we get
\begin{eqnarray*}
A(D_{2})+K A(D_{1})&\geq&\|\varphi\|_{2}^{2}+\|F\|_{2}^{2}\geq
\frac{m(\gamma_{\varphi})+m(\gamma_{F})}{2\pi}\min\{d_{\Omega_{1}}^{2}(0),d_{\Omega_{2}}^{2}(0)\}\\
&=&\min\{d_{\Omega_{1}}^{2}(0),d_{\Omega_{2}}^{2}(0)\}.
\end{eqnarray*}

Next we  prove the sharpness part. We consider the case that
$K=1.$ Let $\Omega_{1}=\Omega_{2}=\mathbb{D}$ and for
$z\in\mathbb{D}$, let $f(z)=tz$, where $0<t<1.$ Then
$$A(D_{1})+A(D_{2})=(1+t^{2})
$$
and $d(0)=1.$ The arbitrariness of $t$
shows the estimate of \eqref{eq-t1} is sharp. The proof of
the theorem is complete.
 \qed

\vspace{6pt}

The following result is well-known.

\begin{Lem}\label{Lem-A}
Among all rectifiable Jordan curves of a given length, the circle
has the maximum interior area.
\end{Lem}

\subsection*{Proof of Theorem \ref{thm-4}} We first prove (\ref{thm-eq1}). By elementary computations, we
have
\beq \label{eqt15} \nonumber
l_{f}(r)&=&r \int_{0}^{2\pi}\left|f_{z}(re^{i\theta})-e^{-2i\theta}
f_{\overline{z}}(re^{i\theta})\right|d\theta\\ \nonumber
&\geq& r\int_{0}^{2\pi}\big(|f_{z}(re^{i\theta})|-|f_{\overline{z}}(re^{i\theta})|\big)d\theta\\
&\geq&\frac{r}{K}\int_{0}^{2\pi}\Lambda_{f}(re^{i\theta})\,d\theta.
\eeq
Cauchy's integral formula applied to $h'(z)=f_z(z)$ and $g'(z)=\overline{f_{\overline{z}}(z)}$
shows for $n\geq1$,
\be\label{eq-t16}
na_{n}=\frac{1}{2\pi i}\int_{|z|=r}\frac{f_{z}(z)}{z^{n}}\,dz~\mbox{ and }~
nb_{n}=\frac{1}{2\pi i}\int_{|z|=r}\frac{\overline{f_{\overline{z}}(z)}}{z^{n}}\,dz,
\ee
respectively.  By (\ref{eqt15}) and (\ref{eq-t16}), we get
\begin{eqnarray*}
n(|a_{n}|+|b_{n}|)&=&\frac{1}{2\pi}\left(\left |\int_{|z|=r}\frac{f_{z}(z)}{z^{n}}\,dz\right |+
\left |\int_{|z|=r}\frac{\overline{f_{\overline{z}}(z)}}{z^{n}}\,dz\right|\right)\\
&\leq&\frac{1}{2\pi
r^{n}}\int_{0}^{2\pi}r\Lambda_{f}(re^{i\theta})\,d\theta\\
&\leq&\frac{Kl_{f}(r)}{2\pi r^{n}}\leq\frac{Kl_{f}(1)}{2\pi r^{n}},
\end{eqnarray*}
which implies that
$$|a_{n}|+|b_{n}|\leq\frac{Kl_{f}(1)}{2n\pi}.
$$

Now we are ready to prove the inequality (\ref{thm-eq2}). First we observe that
\be\label{eq-t9}
S_{f}(r)=\int_{\mathbb{D}_{r}}J_{f}(z)\,d\sigma(z)\geq\frac{1}{K}\int_{\mathbb{D}_{r}}\Lambda_{f}^{2}(z)\,d\sigma(z).
\ee
For $\theta\in[0,2\pi)$ and $z\in\mathbb{D}$, let
$P_{\theta}(z)=\big(f_{z}(z)+e^{i\theta}\overline{f_{\overline{z}}(z)}\big)^{2}$.
By (\ref{eq-t9}) and subharmonicity of $|P_{\theta}|$, we have
\begin{eqnarray*}
|P_{\theta}(z)|&\leq& \frac{1 }{\pi(1-|z|)^{2}}\int_{0}^{1-|z|}\int_{0}^{2\pi}|P_{\theta}(z+\rho
e^{i\beta})|\rho \,d\beta \,d\rho\\
&\leq&\frac{1}{(1-|z|)^{2}}\int_{\mathbb{D}_{1-|z|}}\Lambda_{f}^{2}(z)\,d\sigma(z)\\
&\leq&\frac{S_{f}(1)K}{(1-|z|)^{2}},
\end{eqnarray*}
and  the arbitrariness of $\theta\in[0,2\pi)$ gives the inequality
\be\label{thm-2e}
\Lambda_{f}^{2}(z)\leq\frac{S_{f}(1)K}{(1-|z|)^{2}}.
\ee
By Lemma \Ref{Lem-A}, we get
\be\label{eq-t10}
S_{f}(r)\leq\frac{l_{f}^{2}(r)}{4\pi^{2}}.
\ee
By (\ref{thm-2e}) and (\ref{eq-t10}), we have
$$\Lambda_{f}^{2}(z)\leq\frac{l_{f}^{2}(1)K}{4\pi^{2}(1-|z|)^{2}},
$$
which gives
\be\label{eq-w1}
\Lambda_{f}(z)\leq \frac{l_{f}(1)\sqrt{K}}{2\pi(1-|z|)}.
\ee
Finally, $f\in \mathcal{HB}$ easily follows from (\ref{eq-w1}). The proof of this
theorem is complete. \qed

\section{Bloch  and Lipschitz spaces  on harmonic mappings}\label{csw-sec3}

\begin{lem}\label{lem-t1}
Let $\omega$ be a majorant. For $t>0$, if $\lambda\geq1$, then
\be\label{eq-a1}\omega(\lambda t)\leq\lambda\omega(t).\ee
\end{lem}
\bpf The inequality (\ref{eq-a1}) easily follows from the monotonicity of
$\omega(t)/t$ for $t>0.$ The proof of this lemma is complete.
\epf

\begin{Lem}$($\cite[Lemma 1]{Co}$)$ \label{Lem-t2}
Let $z, ~w$ be complex numbers. Then
$$\max_{\theta\in[0,2\pi]}|w\cos\theta+z\sin\theta|=\frac{1}{2}\big(|w+iz|+|w-iz|\big).
$$
\end{Lem}

\subsection*{Proof of Theorem \ref{thm-CPW}} (a)$\Longleftrightarrow$(c) easily follows
from \cite[Theorem 1.1]{CPMW}. We only need to prove
(a)$\Longleftrightarrow$(b). We first prove that (a)$\Longrightarrow$(b).
Let $z, w\in\mathbb{D}$
with $z\neq w$, and let $\varphi(t)=zt+(1-t)w$, where $t\in[0,1]$. Since
$|\varphi(t)|\leq t|z|+(1-t)|w|$, we see that
\begin{eqnarray*}
1-|\varphi(t)|&=&1-t|z|-|w|+t|w|\\
 &\geq&1-t+|w|(t-1)\\
  &=&(1-t)d(w)
\end{eqnarray*}
and
\begin{eqnarray*}
1-|\varphi(t)|&=&1-t|z|-|w|+t|w|\\
&=&1-t|z|-|w|(1-t)\\
&\geq&1-t|z|-(1-t)\\
&=&td(z).
\end{eqnarray*}
Using the last two inequalities, one has
$$\left(1-|\varphi(t)|\right)^{2}\geq (1-t)t d(w)d(z)
$$
and therefore,  we get
\be\label{eq-t20}
\frac{1}{1-|\varphi(t)|}\leq\frac{1}{\sqrt{(1-t)t d(w)d(z)}}.
\ee

By Lemma \ref{lem-t1} and the inequality (\ref{eq-t20}), for any $z,
w\in\mathbb{D}$ with $z\neq w$, we have
\begin{eqnarray*}
|f(z)-f(w)|&=&\left|\int_{0}^{1}\frac{df}{dt}(\varphi(t))\,dt\right|\quad (\zeta =\varphi (t)=w+t(z-w))\\
&=&\left|(z-w)\int_{0}^{1}f_{\zeta}(\varphi(t))\,dt+(\overline{z}-\overline{w})\int_{0}^{1}f_{\overline{\zeta}}(\varphi(t))\,dt\right|\\
&\leq&|z-w|\int_{0}^{1}\left(|f_{\zeta}(\varphi(t))|+|f_{\overline{\zeta}}(\varphi(t))|\right)dt\\
&\leq&|z-w|\int_{0}^{1}\frac{\Lambda_{f}(\varphi(t))}{\omega\left(\frac{1}{1-|\varphi(t)|}\right)}
\omega\left(\frac{1}{1-|\varphi(t)|}\right)dt
\end{eqnarray*}
\begin{eqnarray*}
&\leq&C|z-w|\int_{0}^{1}\omega\left(\frac{1}{1-|\varphi(t)|}\right)dt\\
&\leq&C|z-w|\int_{0}^{1}\omega\left(\frac{1}{\sqrt{(1-t)t
d(w)d(z)}}\right)dt\\
&\leq&C|z-w|\omega\left(\frac{1}{\sqrt{d(w)d(z)}}\right)\int_{0}^{1}\frac{1}{\sqrt{(1-t)t}}\,dt\\
&=&C|z-w|\omega\left(\frac{1}{\sqrt{d(w)d(z)}}\right)\int_{0}^{\frac{\pi}{2}}\frac{2\sin\theta\cos\theta}
{\sqrt{\sin^{2}\theta\cos^{2}\theta}}\,d\theta\\
&=&C\pi|z-w|\omega\left(\frac{1}{\sqrt{d(w)d(z)}}\right),
\end{eqnarray*}
for some constant $C>0$, which gives
$$\frac{|f(z)-f(w)|}{|z-w|}\leq \pi C\omega\left(\frac{1}{\sqrt{d(w)d(z)}}\right).
$$

Now we prove that (b)$\Longrightarrow$(a). Let
$f=h+\overline{g}$, where $h$ and $g$ are analytic in $\mathbb{D}$.
By Lemma \Ref{Lem-t2}, we obtain
\beq \label{eqt-21}
\nonumber
\max_{\theta\in[0,2\pi]}|f_{x}(z)\cos\theta+f_{y}(z)\sin\theta|&=&\frac{1}{2}
\left(|f_{x}(z)+if_{y}(z)|+|f_{x}(z)-if_{y}(z)|\right)\\ \nonumber
&=&\frac{1}{2}(|2\overline{g_{x}(z)}|+|2h_{x}(z)|)\\ \nonumber
&=&|h'(z)|+|g'(z)|\\
&=&\Lambda_{f}(z).
\eeq
For $r\in(0,1)$ and $\theta\in[0,2\pi]$, let $w=z+re^{i\theta}$.
Then
\beq \label{eqt-22} \nonumber
\lim_{r\rightarrow0+}\left|\frac{f(z)-f(w)}{z-w}\right|&=&\lim_{r\rightarrow0+}\frac{|f(z)-f(z+re^{i\theta})|}{r}\\
\nonumber &=&|f_{x}(z)\cos\theta+f_{y}(z)\sin\theta|\\ \nonumber
&\leq&C\lim_{r\rightarrow0+}\omega\left(\frac{1}{\sqrt{d(z)d(z+re^{i\theta})}}\right)\\
&=&C\omega\left(\frac{1}{d(z)}\right)
\eeq
for some constant $C>0$. By (\ref{eqt-21}) and  (\ref{eqt-22}), we conclude that
$$\Lambda_{f}(z)=\max_{\theta\in[0,2\pi]}|f_{x}(z)\cos\theta+f_{y}(z)\sin\theta|
\leq\max_{\theta\in[0,2\pi]}C\omega\left(\frac{1}{d(z)}\right)=C\omega\left(\frac{1}{d(z)}\right).
$$ Hence (a)$\Longleftrightarrow$(b)$\Longleftrightarrow$(c).
The proof of the theorem is complete. \qed

\subsection*{Proof of Theorem \ref{thm-1}} We first prove the necessity. Since $\Omega$ is a
$L_{\omega}$-extension domain, we see that for any $z,w\in\Omega$,
by using (\ref{eq1.0}), there is a rectifiable curve
$\gamma\subset\Omega$ joining $z$ to $w$ such that
\begin{eqnarray*}
|f(z)-f(w)|&\leq&\int_{\gamma}\Lambda_{f}(\zeta)\,ds(\zeta)
\leq C\int_{\gamma}\frac{\omega\big(d_{\Omega}(\zeta)\big)}{d_{\Omega}(\zeta)}\,ds(\zeta)
\leq C\omega(|z-w|)
\end{eqnarray*}
for some constant $C>0$.

Now we prove the sufficiency. Let $z\in\Omega$ and
$r=d_{\Omega}(z)/2$. 
For all $w\in\mathbb{D}(z,r)$, using (\ref{eq-t1}), we get
$$f(w)=\frac{1}{2\pi}\int_{0}^{2\pi}\mbox{P}(w,re^{i\theta})f(re^{i\theta}+z)\,d\theta,
$$
where
$$\mbox{P}(w,re^{i\theta})=\frac{r^{2}-|w-z|^{2}}{|w-z-re^{i\theta}|^{2}}.
$$
By elementary calculations,  we have
$$\frac{\partial }{\partial
w}\mbox{P}(w,re^{i\theta})=\frac{-(\overline{w}-\overline{z})|w-z-re^{i\theta}|^{2}-(r^{2}-|w-z|^{2})
(\overline{w}-\overline{z}-re^{-i\theta})}{|w-z-re^{i\theta}|^{4}}
$$
and
$$\frac{\partial }{\partial
\overline{w}}\mbox{P}(w,re^{i\theta})=\frac{-(w-z)|w-z-re^{i\theta}|^{2}-(r^{2}-|w-z|^{2})
(w-z-re^{i\theta})}{|w-z-re^{i\theta}|^{4}}.
$$
Then  for all $w\in\mathbb{D}(z,r/2)$,
\begin{eqnarray*}
\left|\frac{\partial }{\partial
w}\mbox{P}(w,re^{i\theta})\right|&\leq&\frac{|w-z||w-z-re^{i\theta}|^{2}+(r^{2}-|w-z|^{2})|w-z-re^{i\theta}|}{|w-z-re^{i\theta}|^{4}}\\
&\leq&\frac{\frac{r}{2}\frac{9r^{2}}{4}+r^{2}\frac{3r}{2}}{\frac{r^{4}}{4}}=\frac{21}{2r}
\end{eqnarray*}
and
$$\left|\frac{\partial }{\partial \overline{w}}\mbox{P}(w,re^{i\theta})\right|\leq\frac{21}{2r},
$$
which implies that
\begin{eqnarray*}
\Lambda_{f}(w)&=&\frac{1}{2\pi}\Big(\Big|\int_{0}^{2\pi}
\frac{\partial }{\partial
w}\mbox{P}(w,re^{i\theta})\big(f(z+re^{i\theta})-f(z)\big)d\theta\Big|\\
&&+\Big|\int_{0}^{2\pi} \frac{\partial }{\partial
\overline{w}}\mbox{P}(w,re^{i\theta})\big(f(z+re^{i\theta})-f(z)\big)d\theta\Big|
\Big)\\
&\leq&\frac{1}{2\pi}\int_{0}^{2\pi}\left(\Big|\frac{\partial
}{\partial w}\mbox{P}(w,re^{i\theta})\Big|+\Big|\frac{\partial
}{\partial
\overline{w}}\mbox{P}(w,re^{i\theta})\Big|\right)\big|f(z+re^{i\theta})-f(z)\big|d\theta\\
&\leq&\frac{21}{2\pi}\int_{0}^{2\pi}\frac{\big|f(z+re^{i\theta})-f(z)\big|}{r}d\theta\\
&\leq&\frac{21C}{2\pi}\frac{\omega(r)}{r}
=\frac{21C}{\pi}\frac{\omega\left(\frac{d_{\Omega}(z)}{2}\right)}{d_{\Omega}(z)}.
\end{eqnarray*}
Since $\omega(t)$ is increasing on $t\in(0,\infty)$, we conclude
that
$$\Lambda_{f}(z)\leq\frac{21C}{\pi}\frac{\omega\left(\frac{d_{\Omega}(z)}{2}\right)}{d_{\Omega}(z)}
\leq\frac{21C}{\pi}\frac{\omega\big(d_{\Omega}(z)\big)}{d_{\Omega}(z)},
$$
for some constant $C>0$. The proof of the theorem is complete.
\qed



\normalsize

\end{document}